\documentclass [11pt]{amsart}
\usepackage{hyperref}
\usepackage {amsmath, amsthm, amssymb, amscd, mathrsfs, graphicx, extarrows, color, url, enumerate}
\usepackage[centering,letterpaper,dvips]{geometry}
\usepackage[all,cmtip]{xy}

\makeatletter
\def\blfootnote{\gdef\@thefnmark{}\@footnotetext}
\makeatother

\newtheorem {theorem}{Theorem}[section]

\newtheorem {conjecture}[theorem]{Conjecture}

\theoremstyle{remark}
\newtheorem {remark}[theorem]{Remark}
\newtheorem {example}[theorem]{Example}

\DeclareFontFamily{U}{mathx}{\hyphenchar\font45}
\DeclareFontShape{U}{mathx}{m}{n}{
      <5> <6> <7> <8> <9> <10>
      <10.95> <12> <14.4> <17.28> <20.74> <24.88>
      mathx10
      }{}
\DeclareSymbolFont{mathx}{U}{mathx}{m}{n}
\DeclareFontSubstitution{U}{mathx}{m}{n}
\DeclareMathAccent{\widecheck}{0}{mathx}{"71}

\def\polhk#1{\setbox0=\hbox{#1}{\ooalign{\hidewidth
    \lower1.5ex\hbox{`}\hidewidth\crcr\unhbox0}}}  
\def\Z {{\mathbb{Z}}}
\def\R {{\mathbb{R}}}
\def\C {{\mathbb{C}}}
\def\Q {{\mathbb{Q}}}

\def\I{\mathfrak{I}}

\def\cp{\mathbb{CP}}

\def\F {\mathbb{F}}
\def\del {\partial}

\def\Sym{\operatorname{Sym}}
\def\Ta {\mathbb{T}_{\alpha}}
\def\Tb {\mathbb{T}_{\beta}}

\DeclareMathOperator{\id}{\operatorname{id}}

\def\spinc{\operatorname{spin}^c}

\def\swf{\operatorname{SWF}}

\def\pin {\operatorname{Pin}(2)}

\def\H {\mathcal{H}}

\newcommand{\s}{\mathfrak{s}}

\def\LE {\mathcal{LE}}
\def\CLE{\mathcal{CLE}}
\def\CW{\mathit{CW}}
\newcommand{\bunderline}[1]{\underline{#1\mkern-2mu}\mkern2mu }

\def\du {\bar{d}}
\def\dl {\bunderline{d}}
\def\deltau {\bar{\delta}}
\def\deltal {\bunderline{\delta}}
\newcommand\alphas{\boldsymbol\alpha}
\newcommand\betas{\boldsymbol\beta}

\def\CF {\mathit{CF}}
\def\HF {\mathit{HF}}
\newcommand\HFhat{\widehat{\HF}}

\newcommand \HFinf {\HF^{\infty}}
\newcommand \CFo {\CF^{\circ}}

\def\HFI {\mathit{HFI}}

\newcommand \HFIo {\HFI^{\circ}}

\def\inv{\iota}
\def\co{\colon\thinspace}
\begin{document}
\title{Homology cobordism and triangulations}

\author{Ciprian Manolescu}


\begin{abstract}
The study of triangulations on manifolds is closely related to understanding the three-dimensional homology cobordism group. We review here what is known about this group, with an emphasis on the  local equivalence methods coming from $\pin$-equivariant Seiberg-Witten Floer spectra and involutive Heegaard Floer homology.\\
\end {abstract}

\maketitle

\section{Triangulations of manifolds}
\blfootnote{Address: Department of Mathematics, University of California, 
Los Angeles, CA 90095}
\blfootnote{Email: cm@math.ucla.edu}
\blfootnote{The author was supported by the NSF grant DMS-1708320.}

A \emph{triangulation} of a topological space $X$ is a homeomorphism $f: |K| \to X$, where $|K|$ is the geometric realization of a simplicial complex $K$. If $X$ is a smooth manifold, we say that the triangulation is smooth if its restriction to every closed simplex of $|K|$ is a smooth embedding. By the work of Cairns \cite{Cairns} and Whitehead \cite{Whitehead}, every smooth manifold admits a smooth triangulation. Furthermore, this triangulation is unique, up to pre-compositions with piecewise linear (PL) homeomorphisms.

The question of classifying triangulations for \emph{topological} manifolds is much more difficult. Research in this direction was inspired by the following two conjectures.

\medskip

\noindent {\bf Hauptvermutung} (Steinitz \cite{Steinitz}, Tietze \cite{Tietze}): {\em Any two triangulations of a space $X$ admit isomorphic subdivisions.}
\medskip

\noindent {\bf Triangulation Conjecture} (based on a remark by Kneser \cite{Kneser}): {\em Any topological manifold admits a triangulation.}
\medskip

Both of these conjectures turned out to be false. The Hauptvermutung was disproved by Milnor \cite{Milnor}, who used Reidemeister torsion to distinguish two triangulations of a space $X$ that is not a manifold. Counterexamples on manifolds came out of the work of Kirby and Siebenmann \cite{KSbook}. (For a nice survey of the mathematics surrounding the Hauptvermutung, see \cite{Ranicki}.) With regard to the triangulation conjecture, counterexamples were shown to exist in dimension $4$ by Casson \cite{Casson}, and in all dimensions $\geq 5$ by the author \cite{beta}.

When studying triangulations on manifolds, a natural condition that one can impose is that the link of every vertex is PL homeomorphic to a sphere. Such triangulations are called {\em combinatorial}, and (up to subdivision) they are equivalent to PL structures on the manifold. 

In dimensions $\leq 3$, every topological manifold admits a unique smooth and a unique PL structure; cf. \cite{Rado, Moise}. In dimensions $\geq 5$, PL structures on topological manifolds were classified in the 1960's. Specifically, building on work of Sullivan \cite{Sullivan} and Casson \cite{Casson0}, Kirby and Siebenmann \cite{KS, KSbook} proved the following:
\begin{itemize}
\item A topological manifold $M$ of dimension $d \geq 5$ admits a PL structure if and only if a certain obstruction class $\Delta(M) \in H^4(M; \Z/2)$ vanishes;
\item For every $d \geq 5$, there exists a $d$-dimensional manifold $M$ such that $\Delta(M)\neq 0$, that is, one without a PL structure;
\item If $\Delta(M)=0$ for a $d$-dimensional manifold $M$ with $d \geq 5$, then the PL structures on $M$ are classified by elements of $H^3(M; \Z/2)$. (This shows the failure of the Hauptvermutung for manifolds.)
\end{itemize}
Finally, in dimension four, PL structures are the same as smooth structures, and the classification of smooth structures is an open problem---although much progress has been made using gauge theory, starting with the work of Donaldson \cite{Donaldson}. Note that Freedman \cite{Freedman} constructed non-smoothable topological four-manifolds, such as the $E_8$-manifold.

We can also ask about arbitrary triangulations of topological manifolds, not necessarily combinatorial. It is not at all obvious that non-combinatorial triangulations of manifolds exist, but they do.

\begin{example}
Start with a triangulation of a non-trivial homology sphere $M^d$ with $\pi_1(M) \neq 1$; such homology spheres exist in dimensions $d \geq 3$. Take two cones on each simplex, to obtain a triangulation of the suspension $\Sigma M$. Repeat the procedure,  to get a triangulation of the double suspension $\Sigma^2 M$. By the double suspension theorem of Edwards \cite{Edwards} and Cannon \cite{Cannon}, the space $\Sigma^2M$ is homeomorphic to $S^{d+2}$. However, the link of one of the final cone points is $\Sigma M$, which is not even a manifold. Thus, $S^{d+2}$ admits a non-combinatorial triangulation. 
\end{example}

\begin{remark}
One can show that any triangulation of a manifold of dimension $\leq 4$ is combinatorial.
\end{remark}

In general, if we triangulate a $d$-dimensional manifold, the link of a $k$-dimensional simplex has the homology of the $(d-k-1)$-dimensional sphere. (However, the link may not be a manifold, as in the example above.) In the 1970's, Galewski-Stern \cite{GS, GS5} and Matumoto \cite{Matumoto} developed the theory of triangulations of high-dimensional manifolds by considering homology cobordism relations between the links of simplices. Their theory involves replacing the links of simplices with PL manifold resolutions, inductively on dimension; cf. \cite{Martin}. In the process we encounter  the $n$-dimensional (piecewise linear) {\em homology cobordism group} $\Theta^n_{\Z}$, 
defined as follows:
$$\Theta^n_{\Z} = \{Y^n \text{ oriented PL manifolds}, H_*(Y) \cong H_*(S^n) \}  /  \sim$$
where the equivalence relation is given by $Y_0 \sim Y_1 \iff$ there exists a compact, oriented, PL  manifold $W^{n+1}$ with $\del W = (-Y_0) \cup Y_1$ and $H_*(W, Y_i; \Z)=0.$ If $Y_0 \sim Y_1$, we say that $Y_0$ and $Y_1$ are {\em homology cobordant}. Summation in $\Theta^n_{\Z}$ is given by connected sum, the standard sphere $S^n$ gives the zero element, and $-[Y]$ is the class of $Y$ with the orientation reversed. This turns $\Theta^n_{\Z}$ into an Abelian group.

It follows from the work of Kervaire \cite{Kervaire} that $\Theta^n_{\Z} = 0$ for $n \neq 3$. On the other hand, the three-dimensional homology cobordism group $\Theta^3_{\Z}$ is nontrivial. To study $\Theta^3_{\Z}$, note that in dimensions three and four, the smooth and PL categories are equivalent. This shows that we can define  $\Theta^3_{\Z}$ in terms of smooth homology spheres and smooth cobordisms.

The easiest way to see that $\Theta^3_{\Z} \neq 0$ is to consider the {\em Rokhlin homomorphism} 
\begin{equation}
\label{eq:rokhlin}
\mu: \Theta^3_\Z \to \Z/2, \ \  \mu(Y) = {\sigma(W)}/{8} \pmod {2},
\end{equation}
where $W$ is any compact, smooth, spin $4$-manifold with boundary $Y$, and $\sigma(W)$ denotes the signature of $W$.  For example, the Poincar\'e sphere $P$ bounds the negative definite plumbing $-E_8$ of signature $-8$, and therefore has $\mu(P) = 1.$ This implies that $P$ is not homology cobordant to $S^3$, and hence $\Theta^3_{\Z} \neq 0$.

Let us introduce the exact sequence:
\begin{equation}
\label{eq:exact}
0\longrightarrow \ker(\mu) \longrightarrow \Theta^3_{\Z} \stackrel{\mu}{\longrightarrow} \mathbb{Z}/2\longrightarrow 0.
\end{equation}

We are now ready to state the results of Galewski-Stern \cite{GS, GS5} and Matumoto \cite{Matumoto} about triangulations of high-dimensional manifolds. They mostly parallel those of Kirby-Siebenmann on combinatorial triangulations. One difference is that, when studying arbitrary triangulations, the natural equivalence relation to consider is {\em concordance}: Two triangulations of the same manifold $M$ are concordant if there exists a triangulation on $M \times [0,1]$ that restricts to the two triangulations on the boundaries $M \times \{0\}$ and $M \times \{1\}$. 

\begin{itemize}

\item A $d$-dimensional manifold $M$ (for $d \geq 5$) is triangulable if and only if $\delta (\Delta(M))=0 \in H^5(M; \ker(\mu)).$ Here, $\Delta(M) \in H^4(M; \Z/2)$ is the Kirby-Siebenmann obstruction to the existence of PL structures, and $\delta: H^4(M; \Z/2) \to H^5(M; \ker(\mu))$ is the Bockstein homomorphism coming from the exact sequence \eqref{eq:exact}. 

\item There exist non-triangulable manifolds in dimensions $\geq 5$ if and only if the exact sequence \eqref{eq:exact} does not split. (In \cite{beta}, the author proved that it does not split.)

\item If they exist, triangulations on a manifold $M$ of dimension $\geq 5$ are classified (up to concordance) by elements in $H^4(M; \ker(\mu))$.
\end{itemize}

The above results provide an impetus for further studying the group $\Theta^3_{\Z}$, together with the Rokhlin homomorphism.

\section{The homology cobordism group}
Since $\Theta^3_{\Z}$ can be defined in terms smooth four-dimensional cobordisms, it is not surprising  that the tools used to better understand it came from gauge theory. Indeed, beyond the existence of the Rokhlin epimorphism, the first progress was made by Fintushel and Stern in \cite{FSorbifolds}, using Yang-Mills theory: 
\begin{theorem}[Fintushel-Stern \cite{FSorbifolds}]
The group $\Theta^3_{\Z}$ is infinite. For example, it contains a $\Z$ subgroup,  generated by the Poincar\'e sphere $\Sigma(2,3,5)$.
\end{theorem}

Their proof involved associating to a Seifert fibered homology sphere $\Sigma(a_1, \dots, a_k)$ a numerical invariant $R(a_1, \dots, a_k)$, the expected dimension of a certain moduli space of self-dual connections. By combining these methods with Taubes' work on end-periodic four-manifolds \cite{TaubesPer}, one obtains a stronger result:
\begin{theorem}[Fintushel-Stern \cite{FSinstanton}, Furuta \cite{FurutaHom}] 
\label{thm:infgen}
The group $\Theta^3_{\Z}$ contains a $\Z^{\infty}$ subgroup. For example, the classes $[\Sigma(2,3,6k-1)], k \geq 1,$ are linearly independent in $\Theta^3_{\Z}$.
\end{theorem}

When $Y$ is a homology three-sphere, the Yang-Mills equations on $\R \times Y$ were used by Floer \cite{FloerInstanton} to construct his celebrated {\em instanton homology}. From the equivariant structure on instanton homology, Fr{\o}yshov \cite{FroyshovYM} defined a homomorphism
$$ h:  \Theta^3_\Z \to \Z,$$
with the property that $h(\Sigma(2,3,5))=1$ (and therefore $h$ is surjective). This implies the following:

\begin{theorem}[Fr{\o}yshov \cite{FroyshovYM}] 
\label{thm:summand}
The group $\Theta^3_\Z$ has a $\Z$ summand, generated by the Poincar\'e sphere $P=\Sigma(2,3,5)$.
\end{theorem}

Since then, further progress on homology cobordism was made using Seiberg-Witten theory and its symplectic-geometric replacement, Heegaard Floer homology. These will be discussed in Sections~\ref{sec:SW} and \ref{sec:HF}, respectively.

In spite of this progress, the following structural questions about $\Theta^3_{\Z}$ remain unanswered:

\medskip
\noindent {\bf Questions:} {\em Does $\Theta^3_{\Z}$ have any torsion? Does it have a $\Z^\infty$ summand? Is it in fact $\Z^\infty$? }
\smallskip

We remark that the existence of a $\Z^{\infty}$ summand could be established by constructing an epimorphism $\Theta^3_{\Z} \to \Z^{\infty}$. In the context of knot concordance, a result of this type was proved by Hom in \cite{Hom}: Using knot Floer homology, she showed the existence of a $\Z^{\infty}$ summand in the smooth knot concordance group generated by topologically slice knots.

\section{Seiberg-Witten theory}
\label{sec:SW}
The Seiberg-Witten equations \cite{SW1, Witten} are a prominent tool for studying smooth four-manifolds. They form a system of nonlinear partial differential equations with a $U(1)$ gauge symmetry; the system is elliptic modulo the gauge action. In dimension three, the information coming from these equations can be packaged into an invariant called {\em Seiberg-Witten Floer homology} (or {\em monopole Floer homology}). This was defined in full generality, for all three-manifolds, by Kronheimer and Mrowka in their book \cite{KMBook}. For rational homology spheres, alternate constructions were given in \cite{MarcolliWang}, \cite{Spectrum}, \cite{FroyshovSW}. Lidman and the author \cite{LMequivalence} showed that the definitions in \cite{Spectrum} and \cite{KMBook} are equivalent.

In many settings, the Seiberg-Witten equations can be used as a replacement for the Yang-Mills equations. For example, from the $S^1$-equivariant structure on Seiberg-Witten Floer homology one can extract an epimorphism 
$$\delta:  \Theta^3_\Z \to \Z,$$
and give a new proof of Theorem~\ref{thm:summand}; see \cite{FroyshovSW, KMBook}. It is not known whether $\delta$ coincides with the invariant $h$ coming from instanton homology. Note that \cite{FroyshovSW, KMBook} use the same notation $h$ for the invariant coming from Seiberg-Witten theory; to prevent confusion with the instanton one, we write $\delta$ here. We use the normalization that $\delta(P)=1$ for the Poincar\'e sphere $P$. 

The construction of Seiberg-Witten Floer homology in \cite{Spectrum} actually gives a refined invariant: an $S^1$-equivariant Floer stable homotopy type, $\swf$, which can be associated to rational homology spheres equipped with $\spinc$ structures.  The definition of $\swf$ was recently generalized to all three-manifolds (in an ``unfolded'' version) by Khandhawit, J. Lin and Sasahira in \cite{KLS}. 

When the $\spinc$ structure comes from a spin structure, the $S^1$ symmetry of the Seiberg-Witten equations (given by constant gauge transformations) can be expanded to a symmetry by the group $\pin$, where 
$$ \pin = S^1 \cup j S^1 \subset \C \oplus j\C = \H.$$
As observed in \cite{beta}, this turns $\swf$ into a $\pin$-equivariant stable homotopy type, and allows us to define a {$\pin$-equivariant Seiberg-Witten Floer homology}. By imitating the construction of the Fr{\o}yshov invariant $\delta$ in this setting, we obtain three new maps
\begin{equation}
\label{eq:abc}
\xymatrixcolsep{1pc}
\xymatrix{
\alpha, \beta, \gamma:  \Theta^3_{\Z}  \ar@{.>}[r] & \Z.}
\end{equation}
These are not homomorphisms (we use the dotted arrow to indicate that), but on the other hand they are related to the Rokhlin homomorphism from \eqref{eq:rokhlin}:
$$\alpha \equiv \beta \equiv \gamma \equiv \mu \ \  (\text{mod } 2).$$

Under orientation reversal, the three invariants behave as follows:
$$\alpha(-Y) = - \gamma(Y), \ \beta(-Y) = -\beta(Y).$$

The properties of $\beta$ suffice to prove the following.

\begin{theorem}[Manolescu \cite{beta}]
\label{thm:beta}
 There are no $2$-torsion elements $[Y] \in \Theta^3_{\Z}$ with $\mu(Y)=1$. Hence, the short exact sequence \eqref{eq:exact} does not split and, as a consequence of \cite{GS, Matumoto}, non-triangulable manifolds exist in every dimension $\geq 5$.
\end{theorem}

Indeed, if $Y$ were a homology sphere with $2[Y]=0  \in \Theta^3_{\Z}$, then $Y$ would be homology cobordant to $- Y$, which would imply that
$$ \beta(Y) = \beta(-Y) = - \beta(Y)\ \Rightarrow \ \beta(Y)=0\ \Rightarrow \ \mu(Y)=0.$$

An alternate construction of $\pin$-equivariant Seiberg-Witten Floer homology, in the spirit of \cite{KMBook} and applicable to all three-manifolds, was given by F. Lin in \cite{FLin}. In particular, this gives an alternate proof of Theorem~\ref{thm:beta}. Lin's theory was further developed in \cite{FLinExact, FLinConnSums, FLinCorrTerms, DaiPin}. 

The invariants $\alpha, \beta, \gamma$ were computed for Seifert fibered spaces by Stoffregen in \cite{Stoffregen} and by F. Lin in \cite{FLinExact}. One application of their calculations is the following result (a proof of which was also announced earlier by Fr{\o}yshov, using instanton homology).

\begin{theorem}[Fr{\o}yshov \cite{FroyshovTalk}, Stoffregen \cite{Stoffregen}, F. Lin \cite{FLinExact}]
There exist homology spheres that are not homology cobordant to any Seifert fibered space.
\end{theorem}

This should be contrasted with a result of Myers \cite{Myers}, which says that every element of $\Theta^3_{\Z}$ can be represented by a hyperbolic three-manifold.

In \cite{Stoffregen2}, Stoffregen studied the behavior of the invariants $\alpha, \beta, \gamma$ under taking connected sums, and used it to give a new proof of the infinite generation of $\Theta^3_{\Z}$. He found a subgroup $\Z^\infty \subset \Theta^3_{\Z}$ generated by the Brieskorn spheres $\Sigma(p, 2p-1, 2p+1)$ for $p \geq 3$ odd. (Compare Theorem~\ref{thm:infgen}.)

In fact, the information in $\alpha, \beta, \gamma, \delta$, and much more, can be obtained from a stronger invariant, a class in the {\em local equivalence group} $\LE$ defined by Stoffregen in \cite{Stoffregen}. To define $\LE$, we first define a {\em space of type SWF} to be a pointed finite $\pin$-CW complex $X$ such that
\begin{itemize}
\item The $S^1$-fixed point set $X^{S^1}$ is $\pin$-homotopy equivalent to $(\tilde \R^s)^+$, where $\tilde \R$ is the one-dimensional representation of $\pin$ on which $S^1$ acts trivially and $j$ acts by $-1$;
\item The action of $\pin$ on $X - X^{S^1}$ is free.
\end{itemize}

The definition is modeled on the properties of the Seiberg-Witten Floer spectra $\swf(Y)$ for homology spheres $Y$. Any $\swf(Y)$ is the formal (de)suspension of a space of type SWF. The condition on the fixed point set comes from the fact that there is a unique reducible solution to the Seiberg-Witten equations on $Y$.

The elements of $\LE$ are equivalence classes $[X]$, where $X$ is a formal (de)suspension of a space of type SWF, and the equivalence relation (called {\em local equivalence}) is given by: $X_1 \sim X_2 \ \iff \ $ there exist $\pin$-equivariant stable maps 
$$ \phi: X_1 \to X_2, \ \ \psi: X_2 \to X_1,$$
which are both $\pin$-equivalences on the $S^1$-fixed point sets. This relation is motivated by the fact that if $Y_1$ and $Y_1$ are homology cobordant, then the induced cobordism maps on Seiberg-Witten Floer spectra give a local equivalence between $\swf(Y_1)$ and $\swf(Y_2)$. 

We can turn $\LE$ into an Abelian group, with addition given by smash product, the inverse given by taking the Spanier-Whitehead dual, and the zero element being $[S^0]$. We obtain a group homomorphism
$$ \Theta^3_{\Z} \to \LE, \ \  [Y] \to [\swf(Y)].$$

The class $[\swf(Y)]\in \LE$ encapsulates all known information from Seiberg-Witten theory that is invariant under homology cobordism. The group $\LE$ is still quite complicated, but there is a simpler version, called the {\em chain local equivalence group} $\CLE$, which involves chain complexes rather than stable homotopy types. The elements of $\CLE$ are modeled on the cellular chain complexes\footnote{When applied to $\swf$, all our chain complexes and homology theories are reduced, but we drop the usual tilde from notation for simplicity.} $C^{\CW}_*(\swf(Y); \F)$ with coefficients in the field $\F=\Z/2$, viewed as modules over 
$$ C^{\CW}_*(\pin; \F)\cong \F[s, j]/ (sj=j^3s, \ s^2=0, j^4=1).$$
and divided by an equivalence relation (called chain local equivalence), similar to the one used in the definition of $\LE$. We have a natural homomorphism
$$ \LE \to \CLE, \ \ [X] \to [C^{\CW}_*(X; \F)].$$

To construct interesting maps from $ \Theta^3_{\Z}$ to $\Z$, one strategy is to factor them through the groups $\LE$ or $\CLE$. Indeed, the Fr{\o}yshov homomorphism $\delta$ can be obtained this way, by passing from chain complexes to the $S^1$-equvariant Borel cohomology, which is a module over 
$$H^*_{S^1}(pt; \F) = H^*(\cp^{\infty}; \F) = \F[U], \ \text{deg}(U)=2.$$

Given the structure of the $S^1$-fixed point set of $\swf(Y)$, one can show that $H^*_{S^1}(\swf(Y); \F)$ is the direct sum of an infinite tower $\F[U]$ and an $\F[U]$-torsion part. The invariant $\delta(Y)$ is set to be $1/2$ the minimal grading in the $\F[U]$ tower. The resulting homomorphism $\delta$ factors as
$$
\xymatrixcolsep{2pc}
\xymatrix{
 \Theta^3_{\Z} \ar[r] &\LE \ar[r] &\CLE \ar[r]^{\delta} & \Z,}$$
Here, by a slight abuse of notation, we also used $\delta$ to denote the final map from $\CLE$ to $\Z$.

The maps $\alpha, \beta, \gamma$ from \eqref{eq:abc} are constructed similarly to $\delta$, but using the $\pin$-equivariant Borel cohomology $H^*_{\pin}(\swf(Y); \F)$. This is a module over
$$H^*_{\pin}(pt; \F) = H^*(B\pin; \F) = \F[q,v]/(q^3), \ \text{deg}(q)=1, \ \text{deg}(v)=4.$$
In this case, if we just consider the $\F[v]$-module structure, we find three infinite towers of the form $\F[v]$, and $\alpha, \beta, \gamma$ are the minimal degrees of elements in this towers, suitably renormalized. We can write
$$
\xymatrixcolsep{3pc}
\xymatrix{
 \Theta^3_{\Z} \ar[r] &\LE \ar[r] &\CLE \ar@{.>}[r]^{\alpha, \beta, \gamma} & \Z.}$$

Two other numerical invariants $
\xymatrixcolsep{1pc}
\xymatrix{
\deltal, \deltau : \CLE \ar@{.>}[r] & \Z}$ can be obtained by considering the $\Z/4$-equivariant Borel cohomology, where $\Z/4$ is the subgroup
$$ \Z/4 = \{1, -1, j, -j\} \subset \pin = \C \oplus j\C.$$

As shown by Stoffregen in \cite{StoffregenNote}, if one considers the Borel homology for other subgroups $G \subset \pin$, one does not get any information beyond that in $\alpha, \beta, \gamma, \delta, \deltal$ and $\deltau$. 

However, one can consider other equivariant generalized cohomology theories. For example, there are invariants $\kappa_i, \ i \in \{0,1\}$ coming from $\pin$-equivariant K-theory (cf. \cite{kg, FurutaLi}), and $\kappa o_i, \ i=0, \dots, 7,$ from $\pin$-equivariant KO-theory \cite{LinKO}. These factor through $\LE$, albeit not through $\CLE$, and have applications to the study of intersection forms of spin four-manifolds with boundary.

In summary, we have a diagram
\begin{equation}
\label{eq:numinv}
\xymatrixcolsep{3pc}
\xymatrix{
 \Theta^3_{\Z} \ar[r] &\mathcal{LE} \ar@{.>}[d]^{\kappa_i, \kappa o_i} \ar[r] &\mathcal{CLE} 
 \ar[r]^{\delta}
 \ar@{.>}[d]^{\deltal,  \deltau}
 \ar@{.>}[dr]^{\alpha, \beta, \gamma} & \Z,\\
  & \Z & \Z & \Z
 }
 \end{equation}

Recall that $\CLE$ was defined using chain complexes with coefficients in $\F=\Z/2$. One could also take coefficients in other fields, say $\Q$ or $\Z/p$ for odd primes $p$. From the corresponding $S^1$-equivariant Borel cohomology (with coefficients in a field of characteristic $p$) one gets homomorphisms
$$ \delta_p: \LE \to \Z.$$
These are different on $\LE$, but it is not known whether they are different when pre-composed with the map $\Theta^3_{\Z} \to \LE$. For every homology sphere for which computations are available, the values of $\delta_p$ are the same for all $p$.

On the other hand, Stoffregen \cite{Stoffregen} showed that the information in chain local equivalence  (for specific Seifert fibered homology spheres) goes beyond that in the numerical invariants from \eqref{eq:numinv}. In fact, using chain local equivalence, he defined an invariant of homology cobordism that takes the form of an Abelian group, called the {\em connected Seiberg-Witten Floer homology}. 

\smallskip
\noindent {\bf Open problem:} {\em Describe the structure of the groups $\LE$ and $\CLE$, and use it to understand more about $\Theta^3_{\Z}$.}
\smallskip

In particular, it would be interesting to construct more homomorphisms from $\LE$ and $\CLE$ to $\Z$, which could perhaps be used to produce new $\Z$ summands in $\Theta^3_{\Z}$. Of special interest is to construct a lift of the Rokhlin homomorphism to $\Z$, as a homomorphism (rather than just as a map of sets, as is the case with $\alpha, \beta, \gamma$). The existence of such a lift would show that $\Theta^3_{\Z}$ has no torsion with $\mu=1$, thus strengthening Theorem~\ref{thm:beta}. In turn, one can show that this would give a simpler criterion for a high-dimensional manifold to be triangulable: the Galewski-Stern-Matumoto class $\delta (\Delta(M)) \in H^5(M; \ker(\mu))$  could be replaced with an equivalent obstruction in $H^5(M; \Z)$.

\section{Heegaard Floer homology and its involutive refinement}
\label{sec:HF}

In a series of papers \cite{HolDisk, HolDiskTwo, HolDiskFour, AbsGraded}, Ozsv\'ath and Szab\'o developed {\em Heegaard Floer homology}: To every three-manifold $Y$ and $\spinc$ structure $\s$, they associated invariants   
$$ \HFhat(Y, \s),\ \HF^+(Y, \s),\ \HF^-(Y, \s),\ \HFinf(Y, \s).$$
These are defined by choosing a pointed Heegaard diagram 
$$ \H = (\Sigma, {\alphas}, {\betas}, z)$$
consisting of the Heegaard surface $\Sigma$ of genus $g$, two sets of attaching curves $\alphas=\{\alpha_1, \dots, \alpha_g\}$, $\betas=\{\beta_1, \dots, \beta_g\}$, and a basepoint $z \in \Sigma$, away from the attaching curves. The attaching curves describe two handlebodies, which put together should give the three-manifold $Y$. One then considers the Lagrangians
$$\Ta = \alpha_1 \times \dots \times \alpha_g, \ \ \Tb=\beta_1\times \dots \times \beta_g$$
inside the symmetric product $\Sym^g(\Sigma)$. The different flavors ($\widehat{\phantom{a}}, +, -, \infty$) of Heegaard Floer homology are versions of the Lagrangian Floer homology $\mathit{HF}(\Ta, \Tb)$.

The construction of Heegaard Floer homology was inspired by Seiberg-Witten theory: the symmetric product is related to moduli spaces of vortices on $\Sigma$. In fact, it has been recently established \cite{KLT1, CGH1, Taubes} that Heegaard Floer homology is isomorphic to the monopole Floer homology from \cite{KMBook}. In view of \cite{LMequivalence}, we obtain a relation to the different  homologies applied to the Seiberg-Witten Floer spectrum $\swf$ (for rational homology spheres). For example, we have
$$\HFhat \cong {H}_*(\swf), \ \HF^+ \cong {H}_*^{S^1}(\swf).$$

Heegaard Floer homology has had numerous applications to low dimensional topology, and is easier to compute than Seiberg-Witten Floer homology. In fact, it was shown to be algorithmically computable; cf. \cite{SarkarWang, LOThat, MOT}.

With regard to homology cobordism, in \cite{AbsGraded} Ozsv\'ath and Szab\'o defined the {\em correction terms} $d(Y, \s)$, which are analogues of the Fr{\o}yshov invariant $\delta$, and give rise to a homomorphism
$$ d: \Theta^3_{\Z} \to \Z.$$
With the usual normalization in Heegaard Floer theory, we have $d = 2\delta$.

One could also ask about recovering $\pin$-equivariant Seiberg-Witten Floer homology and the  invariants $\alpha, \beta, \gamma$ using Heegaard Floer homology. For technical reasons (related to higher order naturality), this seems currently out of reach. However, in \cite{HMinvolutive}, Hendricks and the author developed {\em involutive Heegaard Floer homology}, as an analogue of $\Z/4$-equivariant Seiberg-Witten Floer homology, for the subgroup $\Z/4=\langle j \rangle \subset \pin$. We start by considering the conjugation symmetry on Heegaard Floer complexes $\CFo$ ($\circ \in \{\widehat{\phantom{a}}, +, -, \infty\}$), coming from interchanging the alpha and beta curves, and reversing the orientation of the Heegaard diagram. When $\s$ is self-conjugate (i.e., comes from a spin structure), the conjugation symmetry gives rise to an automorphism
$$\iota: \CFo(Y, \s) \to \CFo(Y, \s),$$
which is a homotopy involution, that is, $\iota^2 \sim \id$. We then define the corresponding involutive Heegaard Floer homology as the homology of the mapping cone of $1+\iota$:
$$ \HFIo(Y, \s) = H_*(\mathit{Cone}(\CFo(Y) \xrightarrow{\phantom{o} (1+\inv) \phantom{o}} \CFo(Y)) ).$$ 
While the usual Heegaard Floer homologies are modules over $H^*_{S^1}(pt) \cong \F[U],$ the involutive versions are modules over $H^*_{\Z/4}(pt) \cong\F[Q, U]/(Q^2)$, with deg$(U)=-2$, \ deg$(Q)=-1$.

\begin{conjecture} 
For every rational homology sphere $Y$ with a self-conjugate $\spinc$ structure $\s$, we have an isomorphism of $\F[Q, U]/(Q^2)$-modules
$$\HFI^+(Y, \s) \cong H_*^{\Z/4}(\swf(Y, \s); \F).$$
\end{conjecture}

From involutive Heegaard Floer homology one can extract invariants $\dl(Y,\s), \du(Y, \s)$, which are the analogues of (twice) the invariants $\deltal, \deltau$ coming from $H^*_{\Z/4}(\swf)$. We get maps
\[
\xymatrixcolsep{1pc}
\xymatrix{
\dl, \du:  \Theta^3_{\Z}  \ar@{.>}[r] & \Z.}
\]

Involutive Heegaard Floer homology has been computed for various classes of three-manifolds, such as large surgeries on alternating knots \cite{HMinvolutive} and the Seifert fibered rational homology spheres $\Sigma(a_1, \dots, a_k$) (or, more generally, almost-rational plumbings) \cite{DaiM}. There is also a connected sum formula for involutive Heegaard Floer homology \cite{HMZ}, and a related connected sum formula for the involutive invariants of knots \cite{ZemkeHFI}. The latter had applications to the study of rational cuspidal curves \cite{BorodzikHom}. 

The calculations of $\dl$ and $\du$ for the above classes of manifolds (and their connected sums) give more constraints on which $3$-manifolds are homology cobordant to each other; see  \cite{HMinvolutive}, \cite{HMZ}, \cite{DaiStoffregen} for several examples. Furthermore, by imitating Stoffregen's arguments from \cite{Stoffregen}, Dai and the author \cite{DaiM} used $\HFI$ to give a new proof that $\Theta^3_\Z$ has a $\Z^\infty$ subgroup.

The chain local equivalence group $\CLE$ admits an analogue in the involutive context, denoted $\I$, whose definition is quite simple. Specifically, we define an {\em $\inv$-complex} to be a pair $(C, \inv)$, consisting of
\begin{itemize}
\item a $\Z$-graded, finitely generated, free chain complex $C$ over the ring $\F[U]$ (deg $U$=-2), such that there is a graded isomorphism 
$U^{-1}H_*(C) \cong \F[U, U^{-1}]$;
\item a grading-preserving chain homomorphism $\iota \co C \to C$, such that $\iota^2 \sim \id$.
\end{itemize}
We say that two $\inv$-complexes $(C, \inv)$ and $(C', \inv')$ are {\em locally equivalent} if there exist (grading-preserving) homomorphisms
$$ F \co C \to C', \ \ G \co C' \to C$$
such that 
$$F \circ \inv \simeq \inv' \circ F,  \ \ \ G \circ \inv' \simeq \inv \circ G,$$
and $F$ and $G$ induce isomorphisms on $U^{-1}H_*$. 

The elements of $\I$ are the local equivalence classes of $\inv$-complexes, and the multiplication in $\I$ is given by
$$ (C, \inv) * (C', \inv') := (C \otimes_{\F[U]} C', \ \inv \otimes \inv' ).$$

As shown in \cite{HMZ}, there is a homomorphism
$$\Theta^3_{\Z} \to \I, \ \ [Y] \to [(\mathit{CF}^-(Y), \iota)],$$
and the maps $d, \dl, \du$ factor through $\mathfrak{I}$. 

Furthermore, in \cite{HHL}, Hendricks, Hom and Lidman extracted from $\I$ a new invariant of homology cobordism, the {\em connected Heegaard Floer homology}, which a summand of Heegaard Floer homology.

\medskip
\noindent {\bf Open problem:} {\em What is $\mathfrak{I}$ as an Abelian group? Can we use it to say more about $\Theta^3_{\Z}$?}

\section{Variations}
So far, we only studied homology cobordisms between integer homology spheres. However, one can define homology cobordisms between two arbitrary three-manifolds $Y_0$ and $Y_1$, by imposing the same conditions on the cobordism, $H_*(W, Y_i; \Z)=0$, $i=0,1.$ Note that, if $Y_0$ is homology cobordant to $Y_1$, then they necessarily have the same homology. The invariants $d, \dl, \du, \alpha, \beta, \gamma$ admit extensions suitable for studying the existence of homology cobordisms between non-homology spheres; see for example \cite[Section 4.2]{AbsGraded}.

We could also weaken the definition of homology cobordism by using homology with coefficients in an Abelian group $A$ different from $\Z$. One gets an {\em $A$-homology cobordism group} $\Theta^3_A$, whose elements are $A$-homology spheres, modulo the relation of $A$-homology cobordism. Observe, for example, that there are natural maps
$$ \Theta^3_{\Z} \to \Theta^3_{\Z/n} \to \Theta^3_{\Q}.$$
Fintushel and Stern \cite{FSexample} showed that the homology sphere $\Sigma(2,3,7)$ bounds a rational ball, whereas it cannot bound an integer homology ball, because $\mu(\Sigma(2,3,7))=1$. This implies that the map $\Theta^3_{\Z} \to  \Theta^3_{\Q}$ is not injective. It is also not surjective, and in fact its cokernel is infinitely generated; cf. \cite{KimLivingston}. In a different direction, Lisca \cite{Lisca} gave a complete description of the subgroup of $\Theta^3_{\Q}$ generated by lens spaces.

One can also construct other versions of homology cobordism by equipping the three-manifolds with $\spinc$ structures, or self-conjugate $\spinc$ structures. Ozsv\'ath and Szab\'o did the former in \cite{AbsGraded}, where they defined a $\spinc$ homology cobordism group $\theta^c$, and showed  that their correction term gives rise to a homomorphism
$$ d: \theta^c \to \Q.$$
The other invariants $\alpha, \beta, \gamma, \dl, \du$ can be similarly extended to maps from a self-conjugate $\spinc$ homology cobordism group (or, more simply, a spin homology cobordism group) to $\Q$. Moreover, on a $\Z/2$-homology sphere there is a unique self-conjugate $\spinc$-structure, which we can use to produce maps
\[
\xymatrixcolsep{1pc}
\xymatrix{
d, \dl, \du, \alpha, \beta, \gamma:  \Theta^3_{\Z/2}  \ar@{.>}[r] & \Q.}
\]

Finally, let us mention that homology cobordism is closely related to knot concordance. Indeed, a  concordance between two knots $K_0, K_1 \subset S^3$ gives rise to a homology cobordism between the surgeries $S^3_m(K_0)$ and $S^3_m(K_1)$, for any integer $m$. It also gives a $\Q$-homology cobordism between the $p^n$-fold cyclic branched covers $\Sigma_{p^n}(K_0)$ and $\Sigma_{p^n}(K_1)$, for any prime $p$ and $n \geq 1$. Thus, one can get knot concordance invariants from homology cobordism invariants, by applying them to surgeries or branched covers. See \cite{AbsGraded, MOwens, Jabuka, HMinvolutive} for examples of this. For a survey of the knot concordance invariants coming from Heegaard Floer homology, we refer to \cite{HomSurvey}.

\medskip
\noindent{\bf Acknowledgements.} I would like to thank Jennifer Hom, Charles Livingston, Robert Lipshitz, Andrew Ranicki, Sucharit Sarkar, Matt Stoffregen and Ian Zemke for comments on a previous version of this paper.

\bibliographystyle{amsalpha}
\bibliography{biblio}

\end{document}